\documentclass[a4paper]{article}

\usepackage[margin=1in]{geometry} 

\usepackage[T1]{fontenc}
\newcommand{\changefont}[3]{
\fontfamily{#1} \fontseries{#2} \fontshape{#3} \selectfont}

\changefont{ptm}{m}{n}

\usepackage{setspace} 
\usepackage{graphicx}

\usepackage{amsfonts}
\usepackage{amsmath}
\usepackage{amssymb}
\usepackage{graphics}
\usepackage{mathrsfs}
\usepackage{color}
\newtheorem{remark}{Remark}[section]

\newtheorem{theorem}{Theorem}[section]

\newtheorem{lemma}{Lemma}[section]

\newtheorem{definition}{Definition}[section]

\long\def\symbolfootnote[#1]#2{\begingroup%
\def\thefootnote{\fnsymbol{footnote}}\footnote[#1]{#2}\endgroup} 

\begin{document}

\begin{center}
\Large \textbf{Continuous-Time and Discrete-Time Quasilinear Systems with Asymptotically Unpredictable Solutions}
\end{center}

\begin{center}
\normalsize \textbf{Mehmet Onur Fen$^{1,}\symbolfootnote[1]{Corresponding Author. E-mail: monur.fen@gmail.com, onur.fen@tedu.edu.tr}$, Fatma Tokmak Fen$^2$} \\
\vspace{0.2cm}
\textit{\textbf{$^1$Department of Mathematics, TED University, 06420 Ankara, Turkey}} \\

\vspace{0.1cm}
\textit{\textbf{$^2$Department of Mathematics, Gazi University, 06560 Ankara, Turkey}} \\
\vspace{0.1cm}
\end{center}

\vspace{0.3cm}

\begin{center}
\textbf{Abstract} 
\end{center}

A novel type of trajectory on semiflows, called asymptotically unpredictable, was proposed by Fen and Tokmak Fen \cite{Fen24}. The presence of sensitivity, which is an indispensable feature of chaotic dynamics, is a crucial property that arises from such trajectories. In the present paper, we show the existence and uniqueness of asymptotically unpredictable solutions for quasilinear systems with delay making benefit of the contraction mapping principle. Additionally, we introduce the notion of an asymptotically unpredictable sequence. It is verified that there exist asymptotically unpredictable sequences which are not unpredictable. Discrete-time equations possessing asymptotically unpredictable orbits are also under investigation. Examples of continuous-time and discrete-time systems with asymptotically unpredictable solutions are provided.  
\vspace{-0.2cm}

\noindent\ignorespaces

\vspace{0.3cm}

\noindent\ignorespaces \textbf{Mathematics Subject Classification:} 34C60, 39A05
 
\vspace{0.2cm}

\noindent\ignorespaces \textbf{Keywords:} Asymptotically unpredictable solutions, systems with delay, discrete equations

\vspace{0.6cm}


\section{Introduction} \label{sec1}

 The investigation of periodic, quasi-periodic, almost periodic, almost automorphic, recurrent, and Poisson stable motions is one of the subjects handled in research studies concerning qualitative theories of differential and discrete equations \cite	{Cheban20,Fink,Kupper02,Laksh,Nemytskii60,Shcherbakov66,Shcherbakov69,Slyusarchuk}. Among these types of motions, the set of Poisson stable ones is the largest in the sense that it comprises all the others. Besides, the notion of Poisson stability can be regarded as a broad generalization of periodicity \cite{Shcherbakov69}. In 2016 a new kind of trajectory, known as unpredictable, was proposed by Akhmet and Fen \cite{Fen16} for semiflows by performing a modification of Poisson stable motions. A crucial feature of such a trajectory is that it leads to Poincar\'e chaos in the associated quasi-minimal set. A distinctness of this chaos type compared to the ones in the sense of Li-Yorke \cite{LiYorke75} and Devaney \cite{Devaney87} is its formation by means of a single trajectory instead of collection of motions. Together with the existence of an uncountable set of everywhere dense Poisson stable trajectories, sensitivity and transitivity are the features of a Poincar\'e chaotic dynamics \cite{Fen16}.

 It was Miller \cite{Miller19} who proposed a generalization of unpredictability to semiflows with arbitrary acting topological monoids. The notion of a topologically unpredictable point for semiflows on topological spaces was proposed by Mahajan et al. \cite{Mahajan24}. It was proved in \cite{Mahajan24} that in compact uniform spaces unpredictability coincides with topological unpredictability. The reader is referred to papers \cite{Thakur20,Thakur21} for results concerning Poincar\'e chaos on product of semiflows. Continuous unpredictable functions, on the other hand, were introduced in study \cite{Fen17a} by means of the Bebutov dynamical system \cite{Sell}. Moreover, unpredictable solutions in differential and discrete equations with applications to neural networks and secure communication were provided in \cite{Fen18, Fen21, Fen23}.

 The behavior of motions in dynamical systems from the asymptotic point of view is important for prediction over a long period of time. The definitions for asymptotically periodic/almost periodic/recurrent/Poisson stable points of dynamical systems were provided in papers \cite{Cheban1977,Shcherbakov77}. Recently, it was found in study \cite{Fen24} that sensitivity can be achieved in a semiflow under a weaker hypothesis compared to unpredictability. This was performed by introducing a new type of trajectory, called asymptotically unpredictable. Such trajectories are necessarily asymptotically Poisson stable and converge to an unpredictable trajectory under the corresponding semiflow \cite{Fen24}.
 
 The main goal of this paper is the investigation of asymptotically unpredictable solutions for non-autonomous systems of delay differential and discrete equations, which are convenient to model various types of real world phenomena \cite{Benczik05,Glass21,Hellwig02,Mohamad03,Pandey24,Zhao16}. In the case of systems with delay, we implement the description of an asymptotically unpredictable function defined on the real axis provided in \cite{Fen24}. It was shown in \cite{Fen24} that a continuous asymptotically unpredictable function can be decomposed as the sum of an unpredictable function and another one converging to zero. Additionally, in the present study we discuss asymptotically unpredictable solutions of discrete-time equations by means of a new definition.

 The paper \cite{Fen18} deals with unpredictable solutions generated by systems of delay differential and discrete equations. The consideration of asymptotically unpredictable solutions instead of unpredictable ones is the main difference of the present study compared to \cite{Fen18}. The existence of asymptotically unpredictable functions and sequences which are not unpredictable manifests the novelty of our results. The main contributions of the present study are as follows.
 \begin{itemize}
 	\item[i.] The existence and uniqueness of asymptotically unpredictable solutions in quasilinear systems of delay differential equations are rigorously proved. The contraction mapping principle is the main tool utilized for that purpose.
 	\item[ii.] We newly introduce asymptotically unpredictable sequences and provide some techniques which are practical to obtain them. 
 	\item[iii.] It is theoretically approved that there exist asymptotically unpredictable sequences which are not unpredictable. In other words, the set of asymptotically unpredictable sequences contain unpredictable ones as a proper subset. This result is instantiated by providing such a sequence.
 	\item[iv.] Asymptotically unpredictable orbits of discrete-time equations are investigated. It is shown under certain assumptions that such equations possess asymptotic unpredictability. Within this framework the discrete analogue of Gronwall's inequality is implemented. 
 \end{itemize}

 \section{Preliminaries} \label{sec2}
 
 Let us denote by $\mathbb T^+$ either the set of non-negative real numbers or the set of non-negative integers. Suppose that $(X,d)$ is a metric space and $\pi: \mathbb T^+ \times X \to X$ is a semiflow on $X$. In other words, the mapping $\pi$ satisfies the following axioms \cite{Sell}: 
 \begin{itemize}
 	\item[i.] $\pi(0,x)=x$ for every $x \in X$,
 	\item[ii.] $\pi(s,x)$ is continuous in the pair of variables $s$ and $x$,
 	\item[iii.] $\pi(s_1,\pi(s_2,x))=\pi(s_1+s_2,x)$ for every $s_1,s_2\in\mathbb T^+$ and for every $x\in X$.
 \end{itemize}
 
  A point $p\in X$ and its trajectory $\{\pi(s,p): s\in\mathbb T^+ \}$ are called unpredictable if there exist a positive number $\varepsilon_0$ and sequences $\{t_n\}_{n\in\mathbb N}$, $\{\tau_n\}_{n\in\mathbb N}$, both of which diverge to infinity, such that $\displaystyle \lim_{n \to \infty} \pi(t_n,p) =p$ and $d[\pi(t_n + \tau_n,p) , \pi(\tau_n,p)] \geq \varepsilon_0$ for every $n \in \mathbb N$ \cite{Fen16}.
 
 The definition of asymptotic unpredictability in semiflows is as follows.
 
 \begin{definition} (\cite{Fen24}) \label{defnasmptoticunp}
 	A point $q \in X$ and its trajectory $\{\pi(s,q): s\in\mathbb T^+\}$ are asymptotically unpredictable if there exists an unpredictable point $p \in X$ such that $\displaystyle \lim_{s \to \infty} d[\pi(s,q), \pi(s,p)] =0$.
 \end{definition}
 
 One of the crucial results demonstrated by Fen and Tokmak Fen \cite{Fen24} is the presence of sensitivity, which means the divergence of initially nearby trajectories, in a dynamics comprising an asymptotically unpredictable point. More precisely, if $q \in X$ is an asymptotically unpredictable point, then there is a point $r$ in the trajectory of $q$ such that the dynamics on the closure of $\{\pi(s,r): s\in\mathbb T^+\}$ is sensitive. It is clear that transitivity is also present in the dynamics since the closure of a trajectory is under discussion. On the other hand, a continuous asymptotically unpredictable function was defined in \cite{Fen24} as an asymptotically unpredictable point of the Bebutov dynamical system, which is identified via the semiflow $\pi:\mathbb R^+ \times \mathcal{C}(\mathbb R) \to \mathcal{C}(\mathbb R)$ satisfying $\pi(t,u)=u_t$, where $u_t(s)=u(t+s)$. Here, $\mathcal{C}(\mathbb R)$ is the set of all continuous functions from $\mathbb R$ into $\mathbb R^m$ equipped with the metric $\rho:\mathcal{C}(\mathbb R) \times \mathcal{C}(\mathbb R) \to \mathbb R^+$ satisfying the equation $$\rho(u,v)=\displaystyle \sum_{j=1}^{\infty}2^{-j} \rho_j(u,v)$$ in which $$\rho_j(u,v)= \min \bigg\lbrace 1,\sup_{-j\leq s\leq j}\left\|u(s)-v(s) \right\| \bigg\rbrace, \ j\in\mathbb N.$$ 
  
 In the remaining parts, $\mathcal{BC}(\mathbb R)$ stands for the set of all uniformly continuous functions defined on $\mathbb R$ with values in $\mathbb R^m$ which are bounded on the whole real axis. Additionally, we make use of the Euclidean norm for vectors and the spectral norm for square matrices.
    
 The subsequent definition of an unpredictable function was given in paper \cite{Fen18}.
    
 \begin{definition} (\cite{Fen18}) \label{defn1}
  A function $\psi \in \mathcal{BC}(\mathbb R)$ is called unpredictable if there exist positive numbers $\varepsilon_0$, $\delta$ and sequences $\{t_n\}_{n\in\mathbb N}$, $\{u_n\}_{n\in\mathbb N}$ both of which diverge to infinity such that $\left\| \psi(t+t_n) - \psi(t)\right\| \to 0$ as $n \to \infty$ uniformly on compact subsets of $\mathbb R$ and $\left\| \psi(t+t_n) - \psi(t)\right\| \geq \varepsilon_0 $ for each $t \in [u_n-\delta, u_n+\delta]$ and $n\in\mathbb N$.	
 \end{definition}
    
  The next definition of an asymptotically unpredictable function in  $\mathcal{BC}(\mathbb R)$ was proposed in study \cite{Fen24}, where it was also deduced via Theorem III.3 of \cite{Sell} that such a function is not periodic. 
    
  \begin{definition} (\cite{Fen24}) \label{defn3}
  A function $\phi \in \mathcal{BC}(\mathbb R)$ is called asymptotically unpredictable if there exist an unpredictable function $\psi \in \mathcal{BC}(\mathbb R)$ and a function $\theta \in  \mathcal{BC}(\mathbb R)$ satisfying $\displaystyle \lim_{t \to \infty} \left\| \theta(t) \right\|=0$ such that $\phi(t) = \psi(t) + \theta(t)$ for every $t\in\mathbb R$.
  \end{definition}
  
  Definition \ref{defn3} is utilized in Section \ref{sec3} for the verification of asymptotic unpredictability in systems with delay. Taking $\theta$ as the zero function in this definition, one can confirm that an unpredictable function defined on the real axis is asymptotically unpredictable. However, the converse is not true. According to the following assertion, which was proved in \cite{Fen24}, there exist asymptotically unpredictable functions that are not unpredictable. 
  
  \begin{lemma} (\cite{Fen24}) \label{lemmaasympto}
  	Suppose that $\psi \in \mathcal{C}(\mathbb R)$ is unpredictable with $\displaystyle \sup_{t \in \mathbb R} \left\| \psi(t) \right\| \leq M$ for some $M \geq 1$. If $\theta \in \mathcal{C}(\mathbb R)$ is a function such that $\displaystyle \lim_{t \to \infty} \left\|\theta(t) \right\|=0 $ and  $\left\|\theta(r_0) \right\|  \geq 4M$ for some $r_0 \in \mathbb R$, then the asymptotically unpredictable function $\phi \in \mathcal{C}(\mathbb R)$ satisfying $\phi(t) = \psi(t) + \theta(t)$, $t \in \mathbb R$, is not unpredictable.
  \end{lemma}
  
Lemma \ref{lemmaasympto} is required in Section \ref{secexamples} for setting up a non-unpredictable asymptotically unpredictable function. The notion of an unpredictable sequence was first proposed by Akhmet and Fen \cite{Fen18} in the following way. 
 
\begin{definition}(\cite{Fen18}) \label{defn2}
	A bounded sequence $\{\psi_i\}_{i\in\mathbb Z}$ in $\mathbb R^p$ is called unpredictable if there exist a positive number $\varepsilon_0$ and sequences $\{\zeta_n\}_{n\in\mathbb N}$, $\{\eta_n\}_{n\in\mathbb N}$ of positive integers both of which diverge to infinity such that $\left\|\psi_{i+\zeta_n} -\psi_i \right\| \to 0$ as $n\to\infty$ for each $i$ in bounded intervals of integers and $\left\|\psi_{\zeta_n + \eta_n} - \psi_{\eta_n}\right\| \geq \varepsilon_0$ for each $n \in \mathbb N$.
\end{definition}

Analogously to Definition \ref{defn3}, we describe an asymptotically unpredictable sequence as follows.

\begin{definition}  \label{defn4}
	A bounded sequence $\{\phi_i\}_{i\in\mathbb Z}$ in $\mathbb R^p$ is called asymptotically unpredictable if there exist an unpredictable sequence $\{\psi_i\}_{i\in\mathbb Z}$ and a sequence $\{\theta_i\}_{i \in \mathbb Z}$ satisfying $\displaystyle \lim_{i \to \infty} \left\| \theta_i \right\|=0$ such that $\phi_i = \psi_i + \theta_i$ for every $i\in\mathbb Z$.
\end{definition}

It is worth noting that the sequence $\{\theta_i\}_{i\in\mathbb Z}$ mentioned in Definition \ref{defn4} is bounded owing to the boundedness of $\{\phi_i\}_{i\in\mathbb Z}$ and $\{\psi_i\}_{i\in\mathbb Z}$. This definition is utilized for the results provided in Sections \ref{sec4} and \ref{sec5}.

\section{Asymptotically Unpredictable Solutions of Retarded Systems} \label{sec3}

In this section, we investigate asymptotically unpredictable solutions of systems with delay in the form
\begin{eqnarray} \label{delaysystem1}
	x'(t) = Ax(t) + f(x(t-\tau)) + \phi(t),	
\end{eqnarray}
where all eigenvalues of the matrix $A \in \mathbb R^m \times \mathbb R^m$ have negative real parts, $f:\mathbb R^m \to \mathbb R^m$ is a continuous function, $\tau$ is a positive number, and $\phi \in \mathcal{BC(\mathbb R)}$.

Because the real parts of eigenvalues of $A$ are negative, there exist numbers $N \geq 1$ and $\lambda >0$ such that $\left\| e^{At}\right\|  \leq N e^{-\lambda t}$ for every $t \geq 0$ \cite{Hale80}.

The following assumptions on system (\ref{delaysystem1}) are required.
\begin{itemize}
	\item[\textbf{(A1)}] There exists a positive number $M_f$ such that $\displaystyle \sup_{x \in \mathbb R^m}\left\| f(x)\right\| \leq M_f$,
	\item[\textbf{(A2)}] There exists a positive number $L_f$ such that $\left\| f(x)-f(\widetilde x)\right\| \leq L_f \left\|x-\widetilde x \right\|$ for all $x,\widetilde x\in\mathbb R^m$,
	\item[\textbf{(A3)}] $\lambda -2 NL_fe^{\lambda \tau/2}>0$.
\end{itemize}

The main result on the existence and uniqueness of an asymptotically unpredictable solution of (\ref{delaysystem1}) is provided in the subsequent theorem.

\begin{theorem} \label{delaymainthm}
	Suppose that the assumptions $(A1)-(A3)$ are fulfilled. If the function $\phi \in \mathcal{BC}(\mathbb R)$ is asymptotically unpredictable, then system (\ref{delaysystem1}) possesses a unique asymptotically unpredictable solution in $\mathcal{BC}(\mathbb R)$. Additionally, this solution is globally exponentially stable.
\end{theorem}

\noindent \textbf{Proof.} Because the function $\phi$ is asymptotically unpredictable, in accordance with Definition \ref{defn1}, there exist an unpredictable function $\psi \in \mathcal{BC}(\mathbb R)$ and a function $\theta \in \mathcal{BC}(\mathbb R)$ with $\left\|\theta(t) \right\| \to 0 $ as $t \to \infty$ such that $\phi(t) = \psi(t) +\theta(t)$ for every $t \in \mathbb R$. Throughout the proof, we denote
$$K_1 = \displaystyle \frac{N^2 (2M_f +M_{\phi} + M_{\psi})}{\lambda - 2 N L_f e^{\lambda \tau /2}}$$
and
$$K_2 = \displaystyle \frac{N}{\lambda - N L_f},$$
where $M_{\phi} = \displaystyle \sup_{t \in \mathbb R} \left\| \phi(t) \right\|$ and $M_{\psi} = \displaystyle \sup_{t \in \mathbb R} \left\| \psi(t) \right\|$. It is worth noting that the numbers $K_1$ and $K_2$ are positive by assumption $(A3)$.

One can confirm using the theory of delay differential equations \cite{Driver1977} that there is a unique globally exponentially stable solution $\Phi \in \mathcal{BC}(\mathbb R)$ of system (\ref{delaysystem1}) which satisfies the relation
\begin{eqnarray} \label{relation1}
	\Phi(t) = \displaystyle \int_{-\infty}^t e^{A(t-s)} \left( f(\Phi(s-\tau)) + \phi(s)\right) ds.
\end{eqnarray}
In a similar way, the retarded system
\begin{eqnarray} \label{delaysystem2}
	y'(t) = Ay(t) + f(y(t-\tau)) + \psi(t),	
\end{eqnarray}
admits a unique globally exponentially stable solution $\Psi \in \mathcal{BC}(\mathbb R)$  such that
\begin{eqnarray} \label{relation2}
	\Psi(t) = \displaystyle \int_{-\infty}^t e^{A(t-s)} \left( f(\Psi(s-\tau)) + \psi(s)\right) ds.
\end{eqnarray}
The equations (\ref{relation1}) and (\ref{relation2}) respectively imply that $\left\|\Phi(t)\right\| \leq \displaystyle \frac{N(M_f + M_{\phi})}{\lambda}$ and $\left\|\Psi(t)\right\| \leq \displaystyle \frac{N(M_f + M_{\psi})}{\lambda}$ for all $t\in\mathbb R$. According to Theorem 2.1 of paper \cite{Fen18}, the solution $\Psi$ of (\ref{delaysystem2}) is unpredictable. In the rest of the proof, we will show that $\left\|\Phi(t) - \Psi(t) \right\| \to 0$ as $t\to \infty$ using the contraction mapping principle.

Let a positive number $\varepsilon$ be given, and suppose that $\gamma$ is a positive number with $\gamma < \displaystyle \frac{1}{K_1 + K_2}$. Owing to the convergence of $\left\|\theta(t) \right\|$ to zero, there is a real number $\alpha$ such that $\left\| \theta(t)\right\|<\gamma \varepsilon$ whenever $t\geq \alpha$.

The function $\Gamma:\mathbb R \to \mathbb R^m$ defined by \begin{eqnarray} \label{funcgamma} 
	\Gamma(t) = \Phi(t) - \Psi(t)
\end{eqnarray}
is a solution of the system
\begin{eqnarray} \label{delaysystem3}
	z'(t) = A z(t) + f(z(t-\tau) + \Psi(t-\tau)) - f(\Psi(t-\tau)) + \theta(t)
\end{eqnarray}
and satisfies the equivalent integral equation
\begin{eqnarray*}
	\Gamma(t) &=& e^{A(t-\alpha)} \left(\Phi(\alpha) - \Psi(\alpha) \right) + \displaystyle \int_{\alpha}^t e^{A(t-s)} \left( f\left(\Gamma(s-\tau) +\Psi(s-\tau) \right)  - f\left(\Psi(s-\tau) \right)  \right) ds \\
	&&+ \displaystyle \int_{\alpha}^t e^{A(t-s)} \theta(s) ds.  
\end{eqnarray*} 

Let $\mathcal{S}$ be the set of all continuous functions $\Gamma$ defined on $\mathbb R$ with values in $\mathbb R^m$ such that $$\left\| \Gamma (t)\right\| \leq K_1 e^{-\lambda (t-\alpha)/2} + K_2 \gamma \varepsilon$$ for $t \geq \alpha -\tau$ and $\left\| \Gamma\right\|_{\infty} \leq M_0$, where
\begin{eqnarray} \label{thenumberm} 
	M_0=\displaystyle \frac{N^2 (2M_f + M_{\phi} + M_{\psi}) + N(M_{\phi} + M_{\psi})}{\lambda - N L_f}
\end{eqnarray}
and $\left\|\Gamma \right\|_{\infty} = \displaystyle \sup_{t \in \mathbb R} \left\| \Gamma (t)\right\|$. On the set $S$, we define an operator $T$ by the equation
$$
T(\Gamma)(t) = \begin{cases}
	\Phi(t) - \Psi(t),  ~ \textrm{ if }  t \leq \alpha, \\
	e^{A(t-\alpha)} \left( \Phi(\alpha) - \Psi(\alpha) \right) + \displaystyle \int_{\alpha}^t e^{A(t-s)} \left( f\left( \Gamma(s-\tau)+\Psi(s-\tau) \right)  - f\left( \Psi(s-\tau) \right) \right) ds\\ + \displaystyle \int_{\alpha}^t e^{A(t-s)} \theta(s)ds, ~  \textrm{ if }  t >\alpha.
\end{cases}
$$

Firstly, we will verify that $T(\mathcal{S}) \subseteq \mathcal{S}$. For that purpose, we fix an element $\Gamma$ of $\mathcal{S}$. For $t > \alpha$, we have 
\begin{eqnarray*}
	\left\|T(\Gamma) (t)\right\| & \leq & \big\|e^{A(t-\alpha)} \big\| \left\|\Phi(\alpha)  - \Psi(\alpha)\right\| + \displaystyle \int_{\alpha}^t \big\|e^{A(t-s)} \big\| \left\| f\left( \Gamma(s-\tau)+\Psi(s-\tau) \right)  - f\left( \Psi(s-\tau) \right)\right\| ds \\
	&&+ \displaystyle \int_{\alpha}^t  \big\|e^{A(t-s)} \big\| \left\|\theta(s)\right\| ds   \\
	&\leq & N e^{-\lambda (t-\alpha)} \left( \displaystyle \frac{N(M_f + M_{\phi})}{\lambda} + \displaystyle \frac{N(M_f + M_{\psi})}{\lambda}\right)  
	+ \displaystyle \int_{\alpha}^t N L_f e^{-\lambda (t-s)} \left\| \Gamma(s-\tau)\right\| ds\\
	&& + \displaystyle \int_{\alpha}^t  N \gamma \varepsilon e^{-\lambda (t-s)} ds \\
	&<&  \displaystyle \frac{N}{\lambda} \left[ N(2M_f +M_{\phi} + M_{\psi} ) + 2 L_f K_1 e^{\lambda \tau/2}  \right] e^{-\lambda (t-\alpha)/2} + \frac{N\gamma \varepsilon}{\lambda}(1+L_f K_2) \\
	&=& K_1 e^{-\lambda(t-\alpha)/2} + K_2 \gamma \varepsilon.
\end{eqnarray*}
Additionally, if $\alpha - \tau \leq t \leq \alpha$, then 
\begin{eqnarray*}
	\left\|T(\Gamma) (t)\right\| = \left\| \Phi(t) - \Psi(t)\right\| \leq \displaystyle \frac{N}{\lambda} (2M_f+M_{\phi} + M_{\psi}) <  K_1.
\end{eqnarray*}
Thus, the inequality $ \left\|T(\Gamma) (t)\right\|\leq K_1 e^{-\lambda(t-\alpha)/2} + K_2 \gamma \varepsilon$ is fulfilled for all $t \geq \alpha-\tau$.

On the other hand, we have for $t > \alpha$ that
\begin{eqnarray*}
	\left\|T(\Gamma) (t)\right\| < \displaystyle \frac{N^2}{\lambda} (2M_f + M_{\phi} + M_{\psi}) + \frac{N}{\lambda} (L_f M_0 + M_{\phi} + M_{\psi}) = M_0,	
\end{eqnarray*}
where $M_0$ is the number defined by (\ref{thenumberm}). Moreover, the inequality 
\begin{eqnarray*}
	\left\|T(\Gamma) (t)\right\| \leq \displaystyle \frac{N}{\lambda} (2M_f + M_{\phi} + M_{\psi}) < M_0,	
\end{eqnarray*}
holds whenever $t \leq \alpha$. Therefore, $\left\| T(\Gamma)\right\|_{\infty} \leq M_0$, and accordingly, $T(\mathcal{S}) \subseteq \mathcal{S}$.

Now, in order to deduce that $T$ is a contraction mapping, let us take two elements $\Gamma$ and $\widetilde{\Gamma}$ of $\mathcal{S}$. One can attain for $t > \alpha$ that
\begin{eqnarray*}
	\big\| T(\Gamma)(t) - T\big(\widetilde{\Gamma}\big)(t)\big\| & \leq & \displaystyle \int_{\alpha}^t N L_f e^{-\lambda (t-s)} \big\| \Gamma(s-\tau)  -  \widetilde{\Gamma} (s-\tau)\big\| ds  
	\leq  \displaystyle \frac{N L_f}{\lambda} \big\| \Gamma - \widetilde{\Gamma}\big\|_{\infty}. 
\end{eqnarray*}
Besides, $T (\Gamma)(t)  - T \big( \widetilde{\Gamma}\big) (t) = 0 $ provided that $ t \leq \alpha$. Thus, $\big\| T (\Gamma)  - T \big( \widetilde{\Gamma}\big) \big\|_{\infty} \leq \displaystyle \frac{N L_f}{\lambda} \big\| \Gamma - \widetilde{\Gamma}\big\|_{\infty}$. The last inequality together with assumption $(A3)$ yield that $T$ is a contraction mapping.

Because of the uniqueness of solutions for the system with delay (\ref{delaysystem3}), the function $\Gamma$ defined by (\ref{funcgamma}) is the unique fixed point of the operator $T$. Let $\Gamma_0$ be the function in $\mathcal{S}$ satisfying
$$
\Gamma_0(t)= \begin{cases}
	\Phi(t) - \Psi(t),  ~ \textrm{if }  t \leq \alpha, \\
	e^{A(t-\alpha)} \left( \Phi(\alpha) - \Psi(\alpha)\right) , ~  \textrm{if }  t >\alpha.
\end{cases}
$$
The sequence of functions $\{\Gamma_k\}$ with $\Gamma_{k+1} = T(\Gamma_k)$, $k=0,1,2,\ldots$, converges to $\Gamma$ on $\mathbb R$. Therefore, 
\begin{eqnarray*}  
	\left\| \Phi(t) - \Psi(t)\right\| \leq K_1 e^{-\lambda(t-\alpha)/2} + K_2 \gamma \varepsilon, \ t \geq \alpha - \tau.
\end{eqnarray*} 
It can be verified for $t \geq \max\Big\lbrace \alpha-\tau, \alpha + \displaystyle \frac{2}{\lambda} \ln \left( \frac{1}{\gamma \varepsilon}  \right) \Big\rbrace $ that  $$\left\| \Phi(t) - \Psi(t)\right\| \leq (K_1 + K_2) \gamma \varepsilon < \varepsilon.$$ Hence, $\left\| \Phi(t) - \Psi(t)\right\| \to 0$ as $t \to \infty$. Consequently, the solution $\Phi$ of (\ref{delaysystem1}) is asymptotically unpredictable. $\square$

\begin{remark} Since there exist asymptotically unpredictable functions in $\mathcal{BC(\mathbb R)}$ which are not unpredictable by Lemma \ref{lemmaasympto}, the result of Theorem \ref{delaymainthm} cannot be reduced to the case of unpredictable solutions. This reveals the novelty of this theorem.
\end{remark}

\section{Results on Asymptotically Unpredictable Sequences} \label{sec4}
 
The content of this section is concerned with the properties of asymptotically unpredictable sequences. Lemma \ref{lemmaseq1} and Lemma \ref{lemmaseq2} provide methods that are convenient to obtain new asymptotically unpredictable sequences from a given one.

\begin{lemma} \label{lemmaseq1}
If $\{\phi_i\}_{i\in\mathbb Z}$ is an asymptotically unpredictable sequence in $\mathbb R^p$, then for every nonsingular matrix $\Omega \in\mathbb R^{p\times p}$ and for every $c \in \mathbb R^p$, the sequence $\{\widetilde{\phi}_i\}_{i\in\mathbb Z}$ satisfying $\widetilde{\phi}_i = \Omega \phi_i+c$, $i\in\mathbb Z$, is also asymptotically unpredictable.
\end{lemma}

\noindent \textbf{Proof.} Suppose that  $\{\psi_i\}_{i\in\mathbb Z}$ is an unpredictable sequence in $\mathbb R^p$ and $\{\theta_i\}_{i\in\mathbb Z}$ is a bounded sequence in $\mathbb R^p$ with $\left\| \theta_i\right\| \to 0$ as $i\to\infty$ such that $\phi_i=\psi_i+\theta_i$ for all $i\in\mathbb Z$. It can be verified that $$\widetilde{\phi}_i=\widetilde{\psi}_i+\widetilde{\theta}_i,\ i\in\mathbb Z,$$ where $\widetilde{\psi}_i=\Omega \psi_i+c$ and $\widetilde{\theta}_i=\Omega \theta_i$.
Due to the unpredictability of $\{\psi_i\}_{i\in\mathbb Z}$, one can find a positive number $\varepsilon_0$ and sequences $\{\zeta_n\}_{n\in\mathbb N}$, $\{\eta_n\}_{n\in\mathbb N}$ of positive integers, both of which diverge to infinity, such that $\left\|\psi_{i+\zeta_n} -\psi_i \right\| \to 0$ as $n \to \infty$ for each $i$ in bounded intervals of integers and $\left\|\psi_{\zeta_n + \eta_n} - \psi_{\eta_n} \right\| \geq \varepsilon_0$, $n \in\mathbb N$.

For each $n\in\mathbb N$ and $i\in\mathbb Z$, we have 
$$\big\|\widetilde{\psi}_{i+\zeta_n} -\widetilde{\psi}_i \big\| \leq \left\| \Omega\right\|\left\|\psi_{i+\zeta_n} -\psi_i \right\| $$
and
$$\big\|\widetilde{\psi}_{\zeta_n + \eta_n} - \widetilde{\psi}_{\eta_n} \big\| \geq \displaystyle \frac{1}{\left\| \Omega^{-1}\right\|} \left\|\psi_{\zeta_n + \eta_n} - \psi_{\eta_n} \right\|.$$
Accordingly, $\big\|\widetilde{\psi}_{i+\zeta_n} -\widetilde{\psi}_i \big\| \to 0$ as $n \to \infty$ for each $i$ in bounded intervals of integers, and $$\big\|\widetilde{\psi}_{\zeta_n + \eta_n} - \widetilde{\psi}_{\eta_n} \big\| \geq \frac{\varepsilon_0}{\left\|\Omega^{-1}\right\|}, \ n \in\mathbb N.$$ Therefore, the sequence $\{\widetilde{\psi}_i\}_{i\in\mathbb Z}$ is unpredictable. Moreover, using the inequality $\big\|\widetilde{\theta}_i\big\| \leq \left\| \Omega\right\| \left\|\theta_i \right\|$, $i\in\mathbb Z$, we obtain that $\displaystyle \lim_{i\to\infty} \big\|\widetilde{\theta}_i \big\| =0$. Thus, the sequence $\{\widetilde{\phi}_i\}_{i\in\mathbb Z}$ is asymptotically unpredictable. $\square$

\begin{lemma}\label{lemmaseq2}
Suppose that $\{\phi_i\}_{i\in\mathbb Z}$ is an asymptotically unpredictable sequence in $\mathbb R^p$ and $\{\varphi_i\}_{i\in\mathbb Z}$ is a bounded sequence in $\mathbb R^p$ such that $\displaystyle \lim_{i\to\infty} \varphi_i$ exists. Then, the sequence $\{\widetilde{\phi}_i\}_{i\in\mathbb Z}$ defined by $\widetilde{\phi}_i=\phi_i+\varphi_i$, $i\in\mathbb Z$, is also asymptotically unpredictable.
\end{lemma}

\noindent \textbf{Proof.} The equation $\phi_i=\psi_i+\theta_i$ is fulfilled for each $i\in\mathbb Z$, where $\{\psi_i\}_{i\in\mathbb Z}$ is an unpredictable sequence in $\mathbb R^p$ and $\{\theta_i\}_{i\in\mathbb Z}$ is a bounded sequence converging to zero. Suppose that $\displaystyle \lim_{i\to\infty} \varphi_i=c$ for some $c\in\mathbb R^p$. Let us define the sequence $\{\widetilde{\theta}_i\}_{i\in\mathbb Z}$ by the equation $\widetilde{\theta}_i=\theta_i+\varphi_i-c$. The inequality $\big\| \widetilde{\theta}_i\big\| \leq \left\| \theta_i\right\| + \left\|\varphi_i-c \right\|$, $i\in\mathbb Z$, yields $\displaystyle \lim_{i\to\infty} \big\| \widetilde{\theta}_i\big\| =0$. Therefore, $\{\psi_i+\widetilde{\theta}_i\}_{i\in\mathbb Z}$ is asymptotically unpredictable. Hence, the same is true for $\{\widetilde{\phi}_i\}_{i\in\mathbb Z}$ by Lemma \ref{lemmaseq1}. $\square$
 
 The preservation of asymptotic unpredictability under shifts is mentioned in the next lemma.
 
\begin{lemma}
If $\{\phi_i\}_{i\in\mathbb Z}$ is an asymptotically unpredictable sequence in $\mathbb R^p$, then for a fixed integer $m$ the sequence $\{\widetilde{\phi}_i\}_{i\in\mathbb Z}$ defined by $\widetilde{\phi}_i=\phi_{i+m}$, $i\in\mathbb Z$, is also asymptotically unpredictable.
\end{lemma}

\noindent \textbf{Proof.} Because the sequence $\{\phi_i\}_{i\in\mathbb Z}$ is asymptotically unpredictable, there exist an unpredictable sequence $\{\psi_i\}_{i\in\mathbb Z}$ and a sequence $\{\theta_i\}_{i\in\mathbb Z}$ satisfying $\displaystyle \lim_{i\to \infty} \left\| \theta_i\right\|=0 $ such that $\phi_i=\psi_i+\theta_i$, $i\in\mathbb Z$. Then, the equation $$\widetilde{\phi}_i=\widetilde{\psi}_i + \widetilde{\theta}_i, \ i\in\mathbb Z,$$ is fulfilled, where $\widetilde{\psi}_i = \psi_{i+m}$ and  $\widetilde{\theta}_i = \theta_{i+m}$. Moreover, $\big\| \widetilde{\theta}_{i}\big\| \to 0$ as $i\to\infty$. In the rest of the proof, it is sufficient to verify that $\{\widetilde{\psi}_i\}_{i\in\mathbb Z}$ is unpredictable.

According to Definition \ref{defn2}, there exist a positive number $\varepsilon_0$ and sequences $\{\zeta_n\}_{n\in\mathbb N}$, $\{\eta_n\}_{n\in\mathbb N}$ of positive integers, both of which diverge to infinity, such that $\left\|\psi_{i+\zeta_n} -\psi_i \right\| \to 0$ as $n\to\infty$ for each $i$ in bounded intervals of integers and $ \left\| \psi_{\zeta_n + \eta_n} - \psi_{\eta_n} \right\|  \geq \varepsilon_0$ for each $n \in \mathbb N$. Let us denote $p_n=\zeta_{n+k}$ and $q_n=\eta_{n+k}-m$, $n\in\mathbb N$, in which $k$ is a fixed sufficiently large natural number with the property that $ \eta_n  >m$ whenever $n > k$. Then, $\{p_n\}_{n\in\mathbb N}$ and $\{q_n\}_{n\in\mathbb N}$ are sequence of positive integers both of which diverge to infinity.

Fix an arbitrary bounded interval of integers $J$ whose least and greatest elements are respectively $l_1$ and $l_2$. Since $\left\|\psi_{i+p_n} -\psi_i \right\| \to 0$ as $n\to\infty$ for each $i$ between $l_1+m$ and $l_2+m$,  we attain that $\big\|\widetilde{\psi}_{i+p_n} -\widetilde{\psi}_i \big\| \to 0$ as $n\to\infty$ for each $i\in J$.
Besides, we have
$$\big\| \widetilde{\psi}_{p_n+q_n} - \widetilde{\psi}_{q_n}\big\| = \big\|  \psi_{p_n+q_n+m} -\psi_{q_n+m}\big\| = \left\| \psi_{\zeta_{n+k}+\eta_{n+k}} - \psi_{\eta_{n+k}}\right\|  \geq \varepsilon_0, \ n\in\mathbb N.$$
This approves the unpredictability of $\{\widetilde{\psi}_i\}_{i\in\mathbb Z}$. 
Consequently,  $\{\widetilde{\phi}_i\}_{i\in\mathbb Z}$ is asymptotically unpredictable. $\square$

We finish Section \ref{sec4} with the following assertion, which reveals that there exist asymptotically unpredictable sequences in $\mathbb R^p$ which are not unpredictable. It is utilized in Section \ref{secexamples} to construct such a sequence. Thus, the set of all asymptotically unpredictable sequences in $\mathbb R^p$ properly includes the set of all unpredictable ones.

\begin{lemma} \label{notunpredict}
Suppose that $\{\psi_i\}_{i\in\mathbb Z}$ is an unpredictable sequence in $\mathbb R^p$ with $\displaystyle \sup_{i\in\mathbb Z} \left\| \psi_i\right\| \leq K$ for some positive number $K$. If $\{\theta_i\}_{i\in\mathbb Z}$ is a bounded sequence in $\mathbb R^p$ such that $\displaystyle \lim_{i\to\infty}\left\| \theta_i\right\| =0$ and $\left\| \theta_m\right\| \geq 4K$ for some integer $m$, then the asymptotically unpredictable sequence $\{\phi_i\}_{i\in\mathbb Z}$ satisfying $\phi_i=\psi_i+\theta_i$, $i\in\mathbb Z$, is not unpredictable.
\end{lemma}

\noindent \textbf{Proof.} Assume that the sequence $\{\phi_i\}_{i\in\mathbb Z}$ is unpredictable. Then, there is a sequence $\{\zeta_n\}_{n\in\mathbb N}$ of positive integers, which diverges to infinity, such that $\left\|\phi_{i+\zeta_n} - \phi_i \right\| \to 0$ as $n \to \infty$ for each $i$ in bounded intervals of integers. 
One can find a natural number $n_0$ such that both of the inequalities
\begin{eqnarray} \label{contradiction}
\left\| \phi_{m+\zeta_{n_0}} - \phi_m\right\| < K
\end{eqnarray} 
and 
$$\left\|\theta_{m+\zeta_{n_0}} \right\| \leq K$$ are valid. On the other hand, we have 
\begin{eqnarray*}
\left\| \phi_{m+\zeta_{n_0}} - \phi_m\right\| \geq \left\| \theta_m \right\| - \left\| \psi_{m + \zeta_{n_0}} \right\|  -  \left\| \psi_m\right\| - \left\| \theta_{m + \zeta_{n_0}}\right\|  \geq K.
\end{eqnarray*}
The last inequality contradicts to (\ref{contradiction}). Therefore, $\{\phi_i\}_{i\in\mathbb Z}$ is not unpredictable. $\square$

Asymptotically unpredictable solutions generated by a class of non-autonomous discrete equations are investigated in the subsequent section.

 \section{Discrete-Time Systems Possessing Asymptotically Unpredictable Orbits} \label{sec5}
 
 This section is devoted to the existence and uniqueness of asymptotically unpredictable orbits of the discrete-time system
 \begin{eqnarray} \label{discretemain}
 w_{i+1} = B w_i + g(w_i) + \phi_i, \ i \in\mathbb Z,
 \end{eqnarray}
 where $B \in \mathbb R^{p \times p}$ is a nonsingular matrix, $g : \mathbb R^p \to \mathbb R^p$ is a continuous function, and $\{\phi_i\}_{i\in\mathbb Z}$ is a bounded sequence.
 
 The following assumptions are needed.
 \begin{itemize}
 	\item[\textbf{(B1)}] There exists a positive number $M_g$ such that $\displaystyle \sup_{w \in \mathbb R^p}\left\| g(w)\right\| \leq M_g$,
 	\item[\textbf{(B2)}] There exists a positive number $L_g$ such that $\left\| g(w)-g(\widetilde w)\right\| \leq L_g \left\|w-\widetilde w \right\|$ for all $w,\widetilde w\in\mathbb R^p$,
 	\item[\textbf{(B3)}] $\left\|B \right\| + L_g <1$.
 \end{itemize}
 
 The main result of  this section is mentioned in the next theorem.
 
 \begin{theorem} \label{mainthm2}
 Suppose that the assumptions $(B1)-(B3)$ hold. Then, the discrete-time system (\ref{discretemain}) possesses a unique asymptotically unpredictable solution provided that the sequence $\{\phi_i\}_{i\in\mathbb Z}$ is asymptotically unpredictable. Moreover, this solution is globally exponentially stable.
 \end{theorem}
 
 \noindent \textbf{Proof.} Owing to the asymptotic unpredictability of $\{\phi_i\}_{i\in\mathbb Z}$, there exist an unpredictable sequence $\{\psi_i\}_{i\in\mathbb Z}$ and a sequence $\{\theta_i\}_{i\in\mathbb Z}$ which converges to zero such that $\phi_i = \psi_i + \theta_i$ for each $i\in\mathbb Z$.
 
 According to the theory of discrete equations \cite{Laksh}, system (\ref{discretemain}) admits a unique globally exponentially stable bounded orbit $\{\Phi_i\}_{i\in\mathbb Z}$ satisfying the relation
 \begin{eqnarray} \label{discphi}
 \Phi_i = \displaystyle \sum_{j=-\infty}^{i} B^{i-j} \left( g(\Phi_{j-1}) + \phi_{j-1}\right). 
 \end{eqnarray}
Similarly, the discrete-time system
  \begin{eqnarray*} \label{discrete2}
 	v_{i+1} = B v_i + g(v_i) + \psi_i, \ i \in\mathbb Z,
 \end{eqnarray*}
 has a unique globally exponentially stable orbit $\{\Psi_i\}_{i\in\mathbb Z}$, which is bounded. Moreover, the equation
 \begin{eqnarray} \label{discpsi}
 	\Psi_i = \displaystyle \sum_{j=-\infty}^{i} B^{i-j} \left( g(\Psi_{j-1}) + \psi_{j-1}\right)
 \end{eqnarray}
 is fulfilled.
The sequence $\{\Psi_i\}_{i\in\mathbb Z}$ is unpredictable according to Theorem 3.1 provided in paper \cite{Fen18},  In the rest of the proof, we will demonstrate that $\left\|\Phi_i-\Psi_i \right\| \to 0$ as $i \to \infty$.

Fix a positive number $\gamma$ such that
$$\gamma \leq \left(\displaystyle \frac{1}{1-\left\|B \right\| -L_g} + \frac{2M_g + M_{\phi} + M_{\psi}}{1-\left\| B\right\| } \right)^{-1},$$ where $M_{\phi} = \displaystyle \sup_{i\in\mathbb Z} \left\|\phi_i \right\|$ and $M_{\psi} = \displaystyle \sup_{i\in\mathbb Z} \left\|\psi_i \right\|$. Let a positive number $\varepsilon$ be given. Since $\left\| \theta_i \right\| \to 0 $ as $i\to\infty$, there is an integer $\alpha$ such that if $i \geq \alpha$, then $\left\|\theta_i \right\|<\gamma \varepsilon$. The equations (\ref{discphi}) and (\ref{discpsi}) yield
\begin{eqnarray*}  
\Phi_i -\Psi_i= \displaystyle \sum_{j=-\infty}^{i} B^{i-j} \left( g(\Phi_{j-1})-g(\Psi_{j-1}) + \theta_{j-1} \right), \ i\in \mathbb Z.
\end{eqnarray*}
Hence, one can attain for $i > \alpha$ that
\begin{eqnarray*}  
\left\| \Phi_i -\Psi_i\right\| & \leq & \displaystyle \sum_{j=-\infty}^{\alpha} \left\| B\right\|^{i-j} \left( 2M_g + M_{\phi} + M_{\psi}\right)  
  + \displaystyle \sum_{j=\alpha+1}^{i} \left\| B\right\|^{i-j} \left( L_g \left\| \Phi_{j-1} -\Psi_{j-1}\right\|+ \gamma \varepsilon\right)  \\
& = & \displaystyle \frac{\gamma \varepsilon}{1-\left\|B \right\|} + (2M_g +M_{\phi} + M_{\psi} - \gamma \varepsilon) \frac{\left\|B\right\|^{i-\alpha}}{1-\left\|B\right\|}  
 + L_g \left\|B\right\|^{i-1} \displaystyle \sum_{j=\alpha}^{i-1} \left\| B\right\|^{-j} \left\|\Phi_j-\Psi_j\right\|.   
\end{eqnarray*}
Accordingly, the inequality
\begin{eqnarray*} 
u_i \leq \displaystyle \frac{\gamma \varepsilon \left\|B \right\|^{-i}}{1-\left\| B\right\|} + (2M_g + M_{\phi} + M_{\psi} - \gamma \varepsilon) \frac{\left\|B \right\|^{-\alpha} }{1-\left\| B\right\|} + \displaystyle \frac{L_g}{\left\| B\right\|} \sum_{j=\alpha}^{i-1} u_j
\end{eqnarray*}
is valid for $i>\alpha$, where $u_i=\left\|B\right\|^{-i} \left\|\Phi_i-\Psi_i\right\|$.
Now, applying the Gronwall type inequality provided in Theorem 4.1.1 \cite{Agarwal00} we obtain that
\begin{eqnarray*}
u_i & \leq & \displaystyle \frac{\gamma \varepsilon \left\|B \right\|^{-i}}{1-\left\| B\right\|} + (2M_g + M_{\phi} + M_{\psi} - \gamma \varepsilon) \frac{\left\|B \right\|^{-\alpha} }{1-\left\| B\right\|} \\
&& + \displaystyle \frac{L_g}{\left\| B\right\|} \sum_{j=\alpha}^{i-1} \left[ \displaystyle \frac{\gamma \varepsilon \left\|B \right\|^{-j}}{1-\left\| B\right\|} + (2M_g + M_{\phi} + M_{\psi} - \gamma \varepsilon) \frac{\left\|B \right\|^{-\alpha} }{1-\left\| B\right\|}\right] \left(1+ \frac{L_g}{\left\| B\right\|} \right)^{i-j-1} \\
&=& \displaystyle \frac{\gamma \varepsilon \left\|B \right\|^{-i}}{1-\left\| B\right\| -L_g} \left[ 1-(\left\| B\right\| +L_g)^{i-\alpha}\right] + (2M_g + M_{\phi} + M_{\psi}) \frac{\left\|B \right\|^{-i} \left(\left\|B \right\|  +L_g\right)^{i-\alpha}}{1-\left\| B\right\| }.
\end{eqnarray*}
Therefore,
\begin{eqnarray*}
\left\| \Phi_i-\Psi_i\right\| \leq \displaystyle \frac{\gamma \varepsilon}{1-\left\| B\right\|-L_g} \left[1-\left(\left\| B\right\| +L_g \right)^{i-\alpha}  \right]  + \frac{2M_g+M_{\phi}+M_{\psi}}{1-\left\|B \right\|} \left( \left\| B\right\|  + L_g \right)^{i-\alpha} 
\end{eqnarray*}
whenever $i>\alpha$. If $i$ is sufficiently large such that
$$i > \max \bigg\lbrace \alpha, \alpha + \ln \left( \frac{1}{\gamma \varepsilon}\right) \left[\ln \left( \frac{1}{\left\|B \right\| +L_g}\right) \right]^{-1}\bigg\rbrace ,$$
then
\begin{eqnarray*}
\left\| \Phi_i-\Psi_i\right\| &<& \displaystyle \frac{\gamma \varepsilon}{1-\left\| B\right\|-L_g}  + \frac{2M_g+M_{\phi}+M_{\psi}}{1-\left\|B \right\|} \left( \left\| B\right\|  + L_g \right)^{i-\alpha} \\
&\leq & \displaystyle \left(  \frac{1}{1-\left\| B\right\|-L_g}  + \frac{2M_g+M_{\phi}+M_{\psi}}{1-\left\|B \right\|}\right) \gamma \varepsilon \\ &\leq& \varepsilon.
\end{eqnarray*}
The last inequality implies that $\left\| \Phi_i-\Psi_i\right\| \to 0$ as $i\to\infty$. Thus,  $\{\Phi_i\}_{i\in\mathbb Z}$ is asymptotically unpredictable. $\square$

\begin{remark}
Because the set of asymptotically unpredictable sequences  includes unpredictable ones as a proper subset by Lemma \ref{notunpredict}, Theorem \ref{mainthm2} cannot be verified by Theorem 3.1 of paper \cite{Fen18}, which is based on unpredictable sequences. 
\end{remark}

In Section \ref{secexamples} we firstly set up a uniformly continuous asymptotically unpredictable function and an asymptotically unpredictable sequence, both of which are not unpredictable. In addition to this, quasilinear systems of delay differential and discrete equations with asymptotically unpredictable solutions are exemplified.

\section{Examples} \label{secexamples}

Let us consider the logistic map 
\begin{eqnarray} \label{logisticmap1}
\kappa_{i+1} = 3.91 \kappa_i(1-\kappa_i), \ i \in\mathbb Z.
\end{eqnarray}
According to Theorem 4.1 proved in study \cite{Fen17}, the map (\ref{logisticmap1}) possesses an unpredictable orbit $\{\kappa^*_i\}_{i\in \mathbb Z}$, which belongs to the unit interval $[0,1]$.  In the examples provided in this section, the sequence $\{\kappa^*_i\}_{i\in \mathbb Z}$ will be utilized as the source of unpredictability. 

In the next subsection we provide an example of an asymptotically unpredictable function defined on the real axis whose construction is different from the one proposed in  \cite{Fen24}. 

\subsection{An Asymptotically Unpredictable Function}

It was demonstrated in Section 4 of \cite{Fen18} that the uniformly continuous function $h:\mathbb R\to\mathbb R$ defined by
\begin{eqnarray} \label{funcht}
h(t)= \displaystyle \int_{-\infty}^t e^{-2(t-s)} \mu(s) ds
\end{eqnarray}
is unpredictable, in which $\mu$ is the piecewise constant function satisfying $\mu(t)=\kappa^*_i$ for $t\in [i,i+1)$, $i\in\mathbb Z$, and $\{\kappa^*_i\}_{i\in \mathbb Z}$ is an unpredictable orbit of (\ref{logisticmap1}) in the interval $[0,1]$. The function $h$ is bounded such that $\displaystyle \sup_{t \in \mathbb R} \left|h(t)\right| \leq 1/2$.

Let us define the function $\psi:\mathbb R \to \mathbb R^2$ as
\begin{eqnarray} \label{funcpsi}
\psi(t) = \begin{pmatrix} 2h(t) \\ h(t)\end{pmatrix}.
\end{eqnarray}
One can confirm using the unpredictability of $h$ that $\psi$ is also unpredictable. Accordingly, the function $\phi:\mathbb R \to \mathbb R^2$ defined by
\begin{eqnarray} \label{aufuncphi}
\phi(t) = \begin{pmatrix} 2h(t) + \displaystyle \frac{3}{1+e^{t}} \\ h(t) - 5 \textrm{sech}(2t) \end{pmatrix}
\end{eqnarray}
is asymptotically unpredictable since $\phi(t)=\psi(t)+\theta(t)$, where
\begin{eqnarray*}  
\theta(t) = \begin{pmatrix} \displaystyle \frac{3}{1+e^{t}} \\  -5 \textrm{sech}(2t) \end{pmatrix}.
\end{eqnarray*}
We would like to point out that  the function $\phi$ is not unpredictable by Lemma \ref{lemmaasympto} since $\displaystyle \sup_{t\in\mathbb R} \left\| \psi(t)\right\| \leq \sqrt 5/2$ and $\left\| \theta(0)\right\| > 2\sqrt 5$.

\subsection{An Asymptotically Unpredictable Sequence}

In this part of the paper, we will give an example of an asymptotically unpredictable sequence, which is not unpredictable. For that purpose, one more time we make benefit of the unpredictable sequence $\{\kappa^*_i\}_{i\in\mathbb Z}$ generated by the map (\ref{logisticmap1}). Let us consider the sequence $\{\psi_i\}_{i\in\mathbb Z}$ in $\mathbb R^2$ satisfying
\begin{eqnarray}
\psi_i = \begin{pmatrix} \kappa^*_i \\  \displaystyle \frac{1}{4} \kappa^*_i \end{pmatrix}, \ i\in\mathbb Z.
\end{eqnarray} 
The unpredictability of $\{\kappa^*_i\}_{i\in\mathbb Z}$ implies the same feature for $\{\psi_i\}_{i\in\mathbb Z}$.
Because $\left\| \theta_i\right\| \to 0$ as $i \to \infty$, where
\begin{eqnarray*}
\theta_i = \begin{pmatrix} \displaystyle \frac{2}{1+i^2} \\ 3e^{-i^2} \end{pmatrix}, \ i\in\mathbb Z,
\end{eqnarray*}
the sequence $\{\phi_i\}_{i\in\mathbb Z}$ with  
\begin{eqnarray} \label{seqphi}
\phi_i = \begin{pmatrix} \kappa^*_i +\displaystyle \frac{2}{1+i^2} \\ \displaystyle \frac{1}{4} \kappa^*_i + 4e^{-i^2} \end{pmatrix} 
\end{eqnarray}
is asymptotically unpredictable. On the other hand, because the inequalities $\displaystyle \sup_{i\in\mathbb Z} \left\|\psi_i \right\| \leq \sqrt{17}/4$ and $\left\|\theta_0 \right\|>\sqrt{17}$ hold, $\{\phi_i\}_{i\in\mathbb Z}$  is not unpredictable in accordance with Lemma \ref{notunpredict}.

The uniformly continuous function $\phi$ defined by (\ref{aufuncphi}) and the sequence $\{\phi_i\}_{i\in\mathbb Z}$  satisfying (\ref{seqphi}) are respectively used in Subsections \ref{subsec1} and \ref{subsec2} to set up systems of delay differential and discrete equations that admit asymptotically unpredictable solutions.

\subsection{A Retarded System with an Asymptotically Unpredictable Trajectory} \label{subsec1}
We consider the system with delay
\begin{eqnarray} \label{exampledelay}
&& x'(t) = x(t) -3y(t) + \displaystyle \frac{1}{6} \arctan(y(t-0.2)) + 2h(t) + \displaystyle \frac{3}{1+e^{t}}, \nonumber \\
&& y'(t) = 5x(t)-5y(t) + \displaystyle \frac{1}{12} \textrm{arccot}(x(t-0.2)) + h(t) - 5 \textrm{sech}(2t),
\end{eqnarray} 
where $t \in\mathbb R$ and $h:\mathbb R \to \mathbb R$ is the unpredictable function defined by equation (\ref{funcht}). 

One can express system (\ref{exampledelay}) in the form of (\ref{delaysystem1}) in which $\tau=0.2$, 
$$A=\begin{pmatrix}  1 & -3 \\ 5 & -5 \end{pmatrix}, \ f(x,y)=\begin{pmatrix}\arctan(y)/6 \\ \textrm{arccot}(x)/12 \end{pmatrix},$$ and
$\phi:\mathbb R \to \mathbb R$ is the asymptotically unpredictable function given by (\ref{aufuncphi}). The eigenvalues of $A$ are $-2\pm i\sqrt 6$ and
\begin{eqnarray} \label{matrixeqn}
e^{A t}= e^{-2t} P \begin{pmatrix} \displaystyle \cos(\sqrt{6} t) & -\sin(\sqrt{6} t) \\ \sin(\sqrt{6} t) & \cos(\sqrt{6} t) \end{pmatrix} P^{-1},
\end{eqnarray}
where $$P=\begin{pmatrix}  0 & 1 \\ -\sqrt{6}/3 & 1 \end{pmatrix}.$$
It can be attained by means of (\ref{matrixeqn}) that $\left\| e^{At}\right\| \leq Ne^{-\lambda t}$ for $t\geq 0$, where $\lambda=2$ and $N=\big\| P \big\| \big\| P^{-1} \big\| =(4+\sqrt{10})/\sqrt{6}$.
The assumptions $(A1)-(A3)$ are fulfilled for (\ref{exampledelay}) with $M_f=\pi\sqrt 2/12$ and $L_f=1/6$. Therefore, there is a unique asymptotically unpredictable solution of (\ref{exampledelay}), which is globally exponentially stable, by Theorem \ref{delaymainthm}.

\subsection{A Discrete-Time System with an Asymptotically Unpredictable Orbit} \label{subsec2}

Let us take into account the discrete system
\begin{eqnarray} \label{exampledisc}
&& x_{i+1} = \displaystyle \frac{1}{4} x_i - \frac{1}{4} y_i + \frac{1}{5} \sin(x_i) + \kappa^*_i +  \frac{2}{1+i^2}, \nonumber \\
&& y_{i+1} =  \displaystyle \frac{1}{2} x_i + \frac{1}{8} y_i + \frac{1}{10} \cos(2y_i) +  \frac{1}{4} \kappa^*_i + 4e^{-i^2}, 
\end{eqnarray} 
where $i\in\mathbb Z$ and $\{\kappa^*_i\}_{i\in\mathbb Z}$ is an unpredictable orbit of the logistic map (\ref{logisticmap1}) belonging to the unit interval $[0,1]$.

System (\ref{exampledisc}) is in the form of (\ref{discretemain}) with
$$B=\begin{pmatrix} 1/4 & -1/4 \\ 1/2 & 1/8 \end{pmatrix}, \ g(x,y)=\begin{pmatrix} \sin x/5 \\ \cos(2y)/10  \end{pmatrix},$$
and $\{\phi_i\}_{i\in\mathbb Z}$ is the asymptotically unpredictable sequence satisfying (\ref{seqphi}). The assumptions $(B1)-(B3)$ are valid for (\ref{exampledisc}) with $M_g=1/\sqrt{20}$, $L_g=1/5$, and $\left\| B \right\| = \sqrt{5}/4$. 
Hence, in accordance with Theorem \ref{mainthm2}, system (\ref{exampledisc}) possesses a unique asymptotically unpredictable solution, and it is globally exponentially stable.

\section{Conclusion} \label{secconclusion}

This study contributes to the qualitative theories of retarded differential as well as discrete equations such that the existence, uniqueness and exponential stability  of asymptotically unpredictable motions in their dynamics are rigorously proved. In both continuous and discrete cases, such trajectories can be decomposed as the sum of an unpredictable motion and another one converging to zero.
 
The concept of strongly unpredictable motions was introduced by Akhmet et al. \cite{Akhmetcarpat}. According to \cite{Akhmetcarpat}, a function $\psi \in \mathcal{BC}(\mathbb R)$ with $\psi=\left(\psi_1,\psi_2,\ldots,\psi_m \right)$ is said to be strongly unpredictable if there exist positive numbers $\varepsilon_0$, $\delta$ and sequences $\{t_n\}_{n\in\mathbb N}$, $\{u_n\}_{n\in\mathbb N}$ both of which diverge to infinity such that $\left\| \psi(t+t_n) - \psi(t)\right\| \to 0$ as $n \to \infty$ uniformly on compact subsets of $\mathbb R$ and $\left\| \psi_i(t+t_n) - \psi_i(t)\right\| \geq \varepsilon_0 $ for each $t \in [u_n-\delta, u_n+\delta]$, $i=1,2,\ldots, m$ and $n\in\mathbb N$. On the other hand, as the discrete analogue of this definition, it was proposed in \cite{Akhmetphys} that a bounded sequence $\{\psi_i\}_{i\in\mathbb Z}$ in $\mathbb R^p$ with $\psi_i=(\psi_i^1,\psi_i^2,\ldots,\psi_i^p)$  is called strongly unpredictable if there exist a positive number $\varepsilon_0$ and sequences $\{\zeta_n\}_{n\in\mathbb N}$, $\{\eta_n\}_{n\in\mathbb N}$ of positive integers both of which diverge to infinity such that $\left\|\psi_{i+\zeta_n} -\psi_i \right\| \to 0$ as $n\to\infty$ for each $i$ in bounded intervals of integers and $\big| \psi^k_{\zeta_n + \eta_n} - \psi^k_{\eta_n}\big| \geq \varepsilon_0$ for each $k=1,2,\ldots,p$ and $n \in \mathbb N$.
 
As a future work, in the light of the following new definitions, one can utilize and develop the techniques provided in Sections \ref{sec3} and \ref{sec5} to investigate the existence and uniqueness of strongly asymptotically unpredictable solutions of differential and discrete equations. 

\begin{definition}  \label{defn5}
	A function $\phi \in \mathcal{BC}(\mathbb R)$ is called strongly asymptotically unpredictable if there exist a strongly unpredictable function $\psi \in \mathcal{BC}(\mathbb R)$ and a function $\theta \in  \mathcal{BC}(\mathbb R)$ satisfying $\displaystyle \lim_{t \to \infty} \left\| \theta(t) \right\|=0$ such that $\phi(t) = \psi(t) + \theta(t)$ for every $t\in\mathbb R$.
\end{definition}

\begin{definition}  \label{defn6}
	A bounded sequence $\{\phi_i\}_{i\in\mathbb Z}$ in $\mathbb R^p$ is called strongly asymptotically  unpredictable if there exist a strongly unpredictable sequence $\{\psi_i\}_{i\in\mathbb Z}$ and a sequence $\{\theta_i\}_{i \in \mathbb Z}$ satisfying $\displaystyle \lim_{i \to \infty} \left\| \theta_i \right\|=0$ such that $\phi_i = \psi_i + \theta_i$ for every $i\in\mathbb Z$.
\end{definition}


\end{document}